\documentclass{article}
\usepackage{graphicx} 
\usepackage[eandd, nonanonymous]{format/style}

\title{On Reduced Rank Regression}
\author{Théo Rudkiewicz$^\dag$, Manon Verbockhaven$^\dag$ \\
    \footnotesize $^\dag$ equal contribution\\
    \footnotesize  firstname.name@inria.fr }
\date{}

\usepackage{natbib}
\usepackage{amsmath}

\usepackage{amsfonts}
\usepackage{amssymb}
\usepackage{caption}
\usepackage{wrapfig}
\usepackage{makecell}
\usepackage{xcolor}
\usepackage{algorithm}
\usepackage{algpseudocode}
\usepackage{bm}
\usepackage{nicefrac}

\newcommand{\overbar}[1]{\mkern1.5mu\overline{\mkern-1.5mu#1\mkern-1.5mu}\mkern1.5mu}

\newcommand{\norm}[1]{\left\lVert#1\right\rVert}

\DeclareMathOperator*{\argmin}{arg\,min \,}

\let\hbar\relax
\newcommand{\hbar}{\overbar{h}}

\DeclareMathOperator{\rank}{rank}

\DeclareMathOperator*{\Tr}{Tr}

\newcommand{\N}{\mathcal{N}}

\newcommand{\R}{\mathbb{R}}

\newcommand{\Shalf}{{\mathbf{S}^{\frac{1}{2}}}}
\newcommand{\Smhalf}{{\mathbf{S}^{-\frac{1}{2}}}}

\newcommand{\myhidetext}[1]{}

\newcommand{\pinv}[1]{#1^{+}}

\newcommand{\trans}[1]{\ensuremath{#1^{\top}}}
\newcommand{\transp}[1]{\trans{\left(#1\right)}}

\newcommand{\trace}[1]{\Tr\left(#1\right)}
\DeclareMathOperator{\Ima}{Im}

\renewcommand{\epsilon}{\varepsilon}
\colorlet{RED}{red}

\newcommand{\scalarp}[1]{\left<#1\right>}

\def\eqref#1{equation~\ref{#1}}

\def\1{\bm{1}}

\def\vepsilon{{\bm{\epsilon}}}

\def\vu{{\bm{u}}}
\def\vv{{\bm{v}}}

\def\mA{{\bm{A}}}
\def\mB{{\bm{B}}}

\def\mD{{\bm{D}}}
\def\mE{{\bm{E}}}
\def\mF{{\bm{F}}}

\def\mI{{\bm{I}}}

\def\mM{{\bm{M}}}
\def\mN{{\bm{N}}}

\def\mP{{\bm{P}}}
\def\mQ{{\bm{Q}}}

\def\mS{{\bm{S}}}

\def\mU{{\bm{U}}}
\def\mV{{\bm{V}}}
\def\mW{{\bm{W}}}
\def\mX{{\bm{X}}}
\def\mY{{\bm{Y}}}

\newcommand{\suchthat}{\text{ s.t. }}

\usepackage[most]{tcolorbox}

\usepackage{booktabs}
\usepackage{multirow}

\usepackage[hidelinks]{hyperref}
\usepackage{cleveref}

\newcommand{\sxyz}{SXYZ}
\newcommand{\projected}{Projected}
\newcommand{\balanced}{Balanced}
\newcommand\thbalanced{
Let $\trans{\mU}\Lambda\mV$ be the reduced SVD of $\Smhalf\mN$ with $r=\rank(\Smhalf\mN)$, $\mU\in\mathbb{R}^{r, p}$, $\mV\in\mathbb{R}^{r, q}$ and $\Lambda\in\mathbb{R}^{r, r}$. Then
    \begin{align*}
        \mA^* &= \Smhalf\transp{\mU_{:k}}\sqrt{\Lambda_{:k, :k}}\\
        \mB^* &= \sqrt{\Lambda_{:k, :k}} \mV_{:k}
    \end{align*}
    is a solution of \Cref{eq:defABstar}. This pair induces the following solution of \Cref{eq:Mstar}:
    \begin{align}
        \mM^*_{\sqrt{\;}} = \Smhalf \left[ \Smhalf\mN \right]_k
    \end{align}
   
}

\newcommand{\thprojected}{
Let $\trans{\mV}\Lambda^2\mV$ be the reduced SVD of $\trans{\mN}\mS^{+}\mN$ with $r=\rank(\trans{\mN}\mS^{+}\mN)$, $\mV\in\mathbb{R}^{r, q}$ and $\Lambda\in\mathbb{R}^{r, r}$. Then
    \begin{align*}
        \mA^* &= \mS^{+}\mN\trans{\mV_{:k}}\\
        \mB^* &= \mV_{:k}
    \end{align*}
    is a solution of \Cref{eq:defABstar}. This pair induces the following solution of \Cref{eq:Mstar}:
    \begin{align}\label{eq:defOLSproj}
        \mM^*_{\perp} &= \mS^{+}\mN\trans{\mV_{:k}}\mV_{:k}
    \end{align}
}

\newcommand{\thSXYZ}{
Let $\trans{\mE}\mD\mF$ be the reduced SVD of $\mX$ with $r=\rank(\mX)$, $\mE\in\mathbb{R}^{r, n}$, $\mF\in \mathbb{R}^{r, p}$ and $\mD\in\mathbb{R}^{r, r}$. Then
\begin{align}
    \mM^*_{\circ} &= \trans{\mF}\mD^{-1}\left[ \mE\mY\right]_k
\end{align}
is a solution of \Cref{eq:Mstar}. With $\trans{\mP}\Gamma\mQ$ the reduced SVD of $\mM_{\circ}^*$, an associated solution of \Cref{eq:defABstar} is $\mA^* = \transp{\mP_{:k}}\sqrt{\Gamma_{:k, :k}}$ and $\mB^* = \sqrt{\Gamma_{:k, :k}}\,\mQ_{:k}$.
}

\newcommand{\propMsqrtequaldMperp}{
Under the common choice of singular vectors of \cref{rem:2},
the matrix $\mM_{\sqrt{}}^*$ of \cref{th:defOLSequi} and the matrix $\mM^*_{\perp}$ of \cref{th:defOLSproj} are equal.
}

\newcommand{\propMperpequalMcirc}{
Under the common choice of singular vectors of \cref{rem:2},
the matrix $\mM_{\perp}^*$ of \cref{th:defOLSproj} and the matrix $\mM^*_{\circ}$ of \cref{th:SXYZ} are equal.
}

\usepackage{amsthm} 

\newtheorem{theorem}{Theorem}
\newtheorem{definition}{Definition}

\newtheorem{proposition}{Proposition}
\newtheorem{lemma}{Lemma}
\newtheorem{remark}{Remark}

\crefname{proof}{proof}{proofs}
\Crefname{proof}{Proof}{Proofs}

\newtheoremstyle{upright}   
  {10pt}                    
  {10pt}                    
  {}                        
  {}                        
  {\bfseries}               
  {.}                       
  { }                       
  {}                        

\theoremstyle{upright}

\tcbuselibrary{theorems}
\definecolor{theoremblue}{HTML}{2196F3} 
\definecolor{theoremred}{HTML}{DE1F10} 
\definecolor{Defblue}{HTML}{0000A0} 
\definecolor{Backgroundblue}{HTML}{0000A0} 
\definecolor{proofpurple}{HTML}{9013FE} 
\definecolor{defyellow}{HTML}{FFEB3B} 
\definecolor{assumptionred}{HTML}{ff3b3b} 
\definecolor{conventionred}{HTML}{ff3b3b} 
\definecolor{remarkblack}{HTML}{000000} 
\definecolor{remarkorange}{HTML}{dc5800} 
\definecolor{remarkgreen}{HTML}{0C7301}
\definecolor{remarkyellow}{HTML}{DEAE10}
\definecolor{examplegreen}{HTML}{00695c} 
\definecolor{schemarouge}{HTML}{EE220C} 
\definecolor{schemableu}{HTML}{4599EF} 

\tcolorboxenvironment{theorem}{
  enhanced jigsaw,
  sharp corners,
  colframe=theoremred,
  colback=theoremred!2,
  coltitle=black,
  fonttitle=\bfseries,
  borderline west={2pt}{0pt}{theoremred},
  before skip=10pt,
  after skip=10pt,
  boxrule=0pt,
  left=10pt,
  right=10pt,
  breakable
}

\tcolorboxenvironment{definition}{
  enhanced jigsaw,
  sharp corners,
  colframe=Defblue,
  colback=Backgroundblue!2,
  coltitle=black,
  fonttitle=\bfseries,
  borderline west={2pt}{0pt}{Defblue},
  before skip=10pt,
  after skip=10pt,
  boxrule=0pt,
  left=10pt,
  right=10pt,
  breakable
}

\tcolorboxenvironment{corollary}{
  enhanced jigsaw,
  sharp corners,
  colframe=theoremblue,
  colback=Backgroundblue!2,
  coltitle=black,
  fonttitle=\bfseries,
  borderline west={2pt}{0pt}{theoremblue},
  before skip=10pt,
  after skip=10pt,
  boxrule=0pt,
  left=10pt,
  right=10pt,
  breakable
}

\tcolorboxenvironment{proposition}{
  enhanced jigsaw,
  sharp corners,
  colframe=theoremblue,
  colback=Backgroundblue!2,
  coltitle=black,
  fonttitle=\bfseries,
  borderline west={2pt}{0pt}{theoremblue},
  before skip=10pt,
  after skip=10pt,
  boxrule=0pt,
  left=10pt,
  right=10pt,
  breakable
}

\tcolorboxenvironment{condition}{
  enhanced jigsaw,
  sharp corners,
  colframe=theoremblue,
  colback=theoremblue!2,
  coltitle=black,
  fonttitle=\bfseries,
  borderline west={2pt}{0pt}{theoremblue},
  before skip=10pt,
  after skip=10pt,
  boxrule=0pt,
  left=10pt,
  right=10pt,
  breakable
}

\tcolorboxenvironment{example}{
  enhanced jigsaw,
  sharp corners,
  colframe=examplegreen,
  colback=examplegreen!2,
  coltitle=black,
  fonttitle=\bfseries,
  borderline west={2pt}{0pt}{examplegreen},
  before skip=10pt,
  after skip=10pt,
  boxrule=0pt,
  left=10pt,
  right=10pt,
  breakable
}

\tcolorboxenvironment{lemma}{
  enhanced jigsaw,
  sharp corners,
  colframe=theoremblue,
  colback=Backgroundblue!2,
  coltitle=black,
  fonttitle=\bfseries,
  borderline west={2pt}{0pt}{theoremblue},
  before skip=10pt,
  after skip=10pt,
  boxrule=0pt,
  left=10pt,
  right=10pt,
  breakable
}

\tcolorboxenvironment{assumption}{
  enhanced jigsaw,
  sharp corners,
  colframe=assumptionred,
  colback=assumptionred!2,
  coltitle=black,
  fonttitle=\bfseries,
  borderline west={2pt}{0pt}{assumptionred},
  before skip=10pt,
  after skip=10pt,
  boxrule=0pt,
  left=10pt,
  right=10pt,
  breakable
}

\tcolorboxenvironment{remark}{
  enhanced jigsaw,
  sharp corners,
  colframe=remarkyellow,
  colback=remarkyellow!2,
  coltitle=black,
  fonttitle=\bfseries,
  borderline west={2pt}{0pt}{remarkyellow},
  before skip=10pt,
  after skip=10pt,
  boxrule=0pt,
  left=10pt,
  right=10pt,
  breakable
}

\newtcolorbox{figurebox}{
  enhanced jigsaw,
  sharp corners,
  colframe=remarkorange,
  colback=remarkorange!2,
  before skip=10pt,
  after skip=10pt,
  boxrule=1pt,
  breakable,
  before skip=10pt,
  after skip=10pt,
}

\newtcolorbox{myproof}{
  enhanced jigsaw,
  sharp corners,
  colframe=schemableu,
  colback=schemableu!0,
  borderline west={1pt}{0pt}{schemableu},
  borderline south={1pt}{0pt}{schemableu},
  before skip=10pt,
  after skip=10pt,
  boxrule=0pt,
  left=10pt,
  right=10pt,
  breakable,
}

\newcounter{proof}
\renewcommand{\theproof}{\arabic{proof}}
\crefname{proof}{proof}{proofs}
\Crefname{proof}{Proof}{Proofs}

\renewenvironment{proof}[1][\proofname]
{
  \refstepcounter{proof} 
  \begin{myproof}
    \textbf{#1~\theproof.}
}
{
  \hfill$\qed$
  \end{myproof}
}

\begin{document}
\maketitle
\begin{abstract}
    We propose a decoupled solution to the Reduced Rank Regression (RRR) problem, i.e.\ a pair of matrices $(\mA^*, \mB^*)$ whose product $\mA^*\mB^*$ minimizes the least squares error under a rank constraint.
    Using our formalism, we rewrite two other solutions from the literature and compare them with ours. In particular, we prove that all these formulations are equivalent. We then analyze their time complexity, and compare their running time and numerical precision on synthetic data, from which we derive practical recommendations. An open-source implementation of the algorithms and experiments is available here: \url{https://gitlab.inria.fr/trudkiew/reduced_rank_regression}.
\end{abstract}

\section{Introduction}

In this article, we are interested in the solution of a rank constrained variant of the linear regression problem which can be formulated as follows.
\begin{definition}
    
Let $\mX \in \mathbb{R}^{n, p}$ and $\mY \in \mathbb{R}^{n, q}$ and $k \in \mathbb{N}$. We define the reduced rank regression (RRR) problem as finding a matrix $\mM^*$ such that:

\begin{equation}\label{eq:Mstar}
    \mM^* \in \argmin_{\mM\in \mathbb{R}^{p, q}}\norm{\mX\mM - \mY} \;\; s.t. \;\; \rank(\mM) \leq k
\end{equation}
where the norm $\norm{.}$ is the Frobenius norm.
Without loss of generality, we consider $k \leq \min(p, q)$.
\end{definition}
By decoupling the variable $\mM$ as $\mM=\mA\mB$, \Cref{eq:Mstar} can be reformulated 
as:
\begin{align}\label{eq:defABstar}
    \mA^*, \mB^* \in \argmin_{\mA, \mB}\norm{\mX\mA\mB - \mY} \;\; s.t. \;\; \mA \in \R^{p,k}, \mB \in \R^{k, q}
\end{align}
Indeed, $\left\{\mA\mB \;|\; \mA \in \R^{p,k}, \mB \in \R^{k, q}\right\} = \left\{\mM\in\R^{p, q} \;|\; \rank(\mM)\leq k\right\}$: the inclusion $\subset$ holds because $\rank(\mA\mB)\leq k$, and the inclusion $\supset$ follows from the reduced SVD of $\mM$, padded with zeros when $\rank(\mM) < k$. Hence, $(\mA^*, \mB^*)$ is a solution of \Cref{eq:defABstar} if and only if $\mA^*\mB^*$ is a solution of \Cref{eq:Mstar}.

\subsection{Notations}
This section is independent of the rest of the paper and presents the notations. 
We denote by $\mI$ the identity matrix. For any integers $k_1, k_2$, we denote by $\mM \in \mathbb{R}^{k_1\times k_2}$ a matrix with $k_1$ rows and $k_2$ columns, and by $\mM^+$ its pseudo-inverse. For $k\leq k_1$ and $k'\leq k_2$, we denote by $\mM_{:k, \;:k'} :=
\left(\mM_{i, j}\right)_{1\leq i \leq k, 1\leq j \leq k'} \in \mathbb{R}^{k\times k'}$ the matrix composed of the $k$ first rows and $k'$ first columns of $\mM$.
We use $\mM_{:k}$, resp. $\mM_{:, :k'}$, for the $k$ first rows, resp. the $k'$ first columns, of $\mM$. For two matrices $\mM_1, \mM_2$ of the same size, $\scalarp{\mM_1, \mM_2} := \trace{\trans{\mM_1}\mM_2}$ denotes their Frobenius inner product, so that $\norm{\mM}^2 = \scalarp{\mM, \mM}$. We call reduced singular value decomposition (SVD) of $\mM$ the matrix product $\trans{\mP}\Gamma\mQ$ with $r=\rank(\mM)$, $\mP\in\mathbb{R}^{r, k_1}$ and $\mQ\in\mathbb{R}^{r, k_2}$ such that $\mP\trans{\mP} = \mQ\trans{\mQ} = \mI_r$, and $\Gamma \in \mathbb{R}^{r, r}$ diagonal with strictly positive entries in non-increasing order. Note that the singular vectors are the \emph{rows} of $\mP$ and $\mQ$. Furthermore, for any positive integer $k$, we denote $\left[\mM\right]_k := \transp{\mP_{:k}}\;\Gamma_{:k, \;:k}\;\mQ_{:k}$, with the convention that $\left[\mM\right]_k =\mM$ for $k\geq r$. 
Finally, if $\mS = \trans{\mQ} \Gamma \mQ$ is symmetric positive semi-definite, we denote by $\mS^{\frac{1}{2}} := \trans{\mQ}\sqrt{\Gamma}\mQ$ its square root and by $\mS^{-\frac{1}{2}} := \trans{\mQ}\pinv{\sqrt{\Gamma}}\mQ = \pinv{\left(\mS^{\frac{1}{2}}\right)}$ the pseudo-inverse of its square root. Note that $\mS^{-\frac{1}{2}}$ is not an inverse when $\mS$ is singular.

\section{Solutions}

In this section, we present a solution for \Cref{eq:defABstar} and compare it
with two other solutions, from \citet{izenmanReducedrankRegressionMultivariate1975} and \citet{xiangOptimalExactLeast2012}, that we have rewritten with our formalism. Although we were unable to find a proof for the solution of \citet{izenmanReducedrankRegressionMultivariate1975}, we still prove that it is equivalent to our balanced solution, as we prove that they induce the same $\mM^*$.

In a first subsection, we present the mathematical formulas, and in a second one, the computational complexity of each solution.
\subsection{Formulas}
Define the (uncentered) second-moment matrices $\mS := \frac{1}{n}\trans{\mX}\mX\in\mathbb{R}^{p, p}$ and $\mN := \frac{1}{n} \trans{\mX}\mY\in\mathbb{R}^{p, q}$. The following theorem holds.
\begin{theorem}[Balanced RRR]\label{th:defOLSequi}
    \thbalanced
    {\color{gray}\hfill \Cref{proof:defOLSequi}}
\end{theorem}

We now rewrite with our formalism the solution proposed by \citet{izenmanReducedrankRegressionMultivariate1975}, for which no proof is available online.
\begin{theorem}[Projected RRR]\label{th:defOLSproj}
    \thprojected
\end{theorem}
\begin{remark}\label{rem:1}
The solution of \Cref{eq:defOLSproj} corresponds to the minimum-norm ordinary least squares solution $\mS^+\mN$, whose rows (vectors of $\R^q$) are projected onto the span of the rows of $\mV_{:k}$.
\end{remark}
\begin{remark}\label{rem:2}
The matrices $\mV$ and $\Lambda$ of \cref{th:defOLSequi} and \cref{th:defOLSproj} can be chosen equal, and so can the matrix of right singular vectors of $\mE\mY$ of \cref{th:SXYZ}, whose singular values are then $\sqrt{n}\Lambda$.\footnote{This follows from the identity $\trans{\mN}\mS^+\mN = \frac{1}{n}\trans{(\mE\mY)}(\mE\mY)$, shown in \cref{prop:circ=Mperpproof}.} In the rest of the paper, we always assume such a common choice.
\end{remark}
Although we did not find any proof of \Cref{th:defOLSproj}, it follows from \Cref{th:defOLSequi} and the following proposition, since the objective of \Cref{eq:defABstar} only depends on the product $\mA\mB$.
\begin{proposition}\label{prop:Msqrt=Mperp}
    \propMsqrtequaldMperp
    {\color{gray}\hfill \Cref{prop:Msqrt=Mperpproof}}
\end{proposition}

\vspace{0.3cm}
Finally, \citet{xiangOptimalExactLeast2012} provide the following solution of \Cref{eq:Mstar}.
\begin{theorem}[SXYZ RRR]\label{th:SXYZ}
    \thSXYZ
    {\color{gray}\hfill \Cref{proof:SXYZ}}
\end{theorem}
We conclude this section with the following proposition.
\begin{proposition}
    \propMperpequalMcirc
    {\color{gray}\hfill \Cref{prop:circ=Mperpproof}}
\end{proposition}



\subsection{Algorithms and complexities}\label{sec:algorithms}

This sub-section presents the algorithms that compute the three solutions of \Cref{th:defOLSequi,th:defOLSproj,th:SXYZ} and analyzes their time complexity. The main conclusion, summarized in \Cref{tab:complexity}, is that \projected{} and \balanced{} scale as $np^2$ for large $n$ while \sxyz{} scales as $np^2$ only when $p \gg n$. In other words, \sxyz{} is expected to be the fastest method when $n < p$ and \projected{} when $n \geq p$. The experiments of \Cref{sec:experiments} confirm this prediction.

Throughout this sub-section, complexities are given in floating-point operations (flops) and constant factors are dropped, except for the constants $C_{\text{svd}}$ and $C_{\text{eigh}}$ of the decompositions, which are large enough to matter in practice.

\subsubsection*{Cost of elementary operations}

The algorithms rely on four linear algebra operations:
\begin{itemize}
    \item \textbf{Matrix multiplication} of matrices of size $(m \times k)$ and $(k \times n)$: $2mkn$ flops, counted as $mkn$ below since we drop the constant $2$;
    \item \textbf{Singular value decomposition} (\texttt{torch.linalg.svd}) of a matrix of size $(m \times n)$: \\$C_{\text{svd}}\, m\, n\, \min(m,n)$ flops, where $C_{\text{svd}}$ is a constant of order $10$;
    \item \textbf{Symmetric eigendecomposition} (\texttt{torch.linalg.eigh}) of a matrix of size $(m \times m)$: $C_{\text{eigh}}\, m^3$ flops, where $C_{\text{eigh}}$ is also a constant of order $10$;
    \item \textbf{Pseudo-inverse} (\texttt{torch.linalg.pinv}): computed through an SVD, or through \texttt{eigh} when called with \texttt{hermitian=True}, and therefore with the same cost as the corresponding decomposition.
\end{itemize}

\subsubsection*{Common first step: computation of $\trans{\mX}\mX$ and $\trans{\mX}\mY$}

\projected{} and \balanced{} both start by computing the second-moment matrices $\mS \propto \trans{\mX}\mX$ and $\mN \propto \trans{\mX}\mY$. This step costs $n(p^2 + pq)$ flops and is the only one whose cost depends on $n$. It can be done by batches of $b$ samples, which reduces the memory footprint from $O(n(p + q) + p^2 + pq)$ for the naive approach to $O(b(p + q) + p^2 + pq)$. This is why these two methods can handle datasets that do not fit in memory, unlike \sxyz{} which needs an SVD of the full matrix $\mX$.

The algorithms below omit the factor $\frac{1}{n}$ in the definition of $\mS$ and $\mN$, which the implementation keeps. This does not change the product $\mM = \mA\mB$.

\subsubsection*{Projected RRR (\Cref{th:defOLSproj})}

\begin{algorithm}[ht]
\caption{Projected RRR}\label{algo:projected}
\begin{algorithmic}[1]
\Function{Projected-RRR}{$\mX \in \mathbb{R}^{n \times p}$, $\mY \in \mathbb{R}^{n \times q}$, $k$}
    \State \Comment{Compute sufficient statistics \hfill $n(p^2 + pq)$}
    \State $\mS \gets \trans{\mX} \mX$
    \State $\mN \gets \trans{\mX} \mY$
    \State
    \State $\mS^+ \gets \texttt{torch.linalg.pinv}(\mS, \texttt{hermitian=True})$ \Comment{$C_{\text{eigh}}\, p^3$}
    \State $\mP \gets \trans{\mN} \mS^+ \mN$ \Comment{$qp^2 + q^2 p$}
    \State $s, \mV \gets \texttt{eigh}(\mP)$ \Comment{$C_{\text{eigh}}\, q^3$}
    \State \Comment{Keep the $k$ largest eigenvectors}
    \State $\mV \gets \mV[:, -k:]$
    \State $\mA \gets \mS^+ \mN \mV$ \Comment{$p^2 q + pqk$}
    \State $\mB \gets \trans{\mV}$
    \State \Return $\mA, \mB$
\EndFunction
\end{algorithmic}
\end{algorithm}

\Cref{algo:projected} does the low rank selection in the output space: after the second-moment matrices, it computes the pseudo-inverse of $\mS$ (a $p \times p$ eigendecomposition), forms the $q \times q$ matrix $\mP = \trans{\mN}\mS^+\mN$ and keeps its $k$ leading eigenvectors. The factors $\mA$ and $\mB$  are obtained after few matrix products. Adding up the costs of the lines of \Cref{algo:projected} gives a total complexity of
$$ n(p^2 + pq) +  C_{\text{eigh}}(p^3 + q^3) + 2qp^2 + q^2 p + pqk. $$
For $p \approx q \approx k$, this simplifies to $2np^2 + (2 C_{\text{eigh}} + 4)\, p^3 \approx np^2 + 20p^3$: a term linear in $n$ from the second-moment matrices, plus a term independent of $n$ dominated by the two eigendecompositions.

\subsubsection*{Balanced RRR (\Cref{th:defOLSequi})}

\begin{algorithm}[ht]
\caption{Balanced RRR}\label{algo:balanced}
\begin{algorithmic}[1]
\Function{Balanced-RRR}{$\mX \in \mathbb{R}^{n \times p}$, $\mY \in \mathbb{R}^{n \times q}$, $k$}
    \State \Comment{Compute sufficient statistics \hfill $n(p^2 + pq)$}
    \State $\mS \gets \trans{\mX} \mX$
    \State $\mN \gets \trans{\mX} \mY$
    \State
    \State $\Smhalf \gets \texttt{sqrt\_inverse\_matrix}(\mS)$ \Comment{$C_{\text{eigh}}\, p^3$}
    \State $\mP \gets \Smhalf\, \mN$ \Comment{$p^2 q$}
    \State $\mU, \bm{\sigma}, \mV \gets \texttt{torch.linalg.svd}(\mP)$ \Comment{$C_{\text{svd}}\, pq\min(p,q)$}
    \State \Comment{Truncate to rank $k$}
    \State $\mU, \bm{\sigma}, \mV \gets \mU_{:,:k},\; \bm{\sigma}_{:k},\; \mV_{:k,:}$
    \State
    \State $\mA \gets \Smhalf\, \mU \sqrt{\bm{\sigma}}$ \Comment{$p^2 k$}
    \State $\mB \gets \sqrt{\bm{\sigma}}\, \mV$
    \State \Return $\mA, \mB$
\EndFunction
\end{algorithmic}
\end{algorithm}

\Cref{algo:balanced} shares the first step with \Cref{algo:projected} but then follows \Cref{th:defOLSequi} directly: it computes the inverse square root $\Smhalf$ (a $p \times p$ eigendecomposition), whitens $\mN$ into $\mP = \Smhalf\mN$ and takes the rank-$k$ truncated SVD of $\mP$. The total complexity is
$$ n(p^2 + pq) + C_{\text{eigh}} p^3 + C_{\text{svd}} pq\min(p,q) + p^2(k+q), $$
which for $p \approx q \approx k$ simplifies to $2np^2 + (C_{\text{eigh}} + C_{\text{svd}} + 2) p^3 \approx np^2 + 20p^3$. The asymptotic cost is therefore the same as \projected{}, the second eigendecomposition being replaced by an SVD of the same size. In practice, because the SVD 
is more expensive than the eigendecomposition, we experimentally see in \Cref{sec:experiments} that
\balanced{} is slower than \projected{}.

\subsubsection*{SXYZ RRR (\Cref{th:SXYZ})}

\begin{algorithm}[ht]
\caption{SXYZ RRR}\label{algo:sxyz}
\begin{algorithmic}[1]
\Function{SXYZ-RRR}{$\mX \in \mathbb{R}^{n \times p}$, $\mY \in \mathbb{R}^{n \times q}$, $k$}
    \State $\mU, \bm{\sigma}, \mV \gets \texttt{torch.linalg.svd}(\mX)$ \Comment{$C_{\text{svd}}\, np\min(n,p)$}
    \State \Comment{Keep the $r$ singular values above the cut-off}
    \State $\mU, \bm{\sigma}, \mV \gets \mU_{:,:r},\; \bm{\sigma}_{:r},\; \mV_{:r,:}$
    \State $\mW' \gets \trans{\mU} \mY$ \Comment{$nq\min(n, p)$}
    \State $\mW \gets \texttt{svd\_projection}(\mW', k)$ \Comment{$C_{\text{svd}}\, q \min(n, p, q) \min(n,p)$}
    \State $\mM \gets \trans{\mV} (\bm{\sigma}^{-1} \odot \mW)$ \Comment{$pq\min(n, p)$}
    \State \Return $\mM$
\EndFunction
\end{algorithmic}
\end{algorithm}

Unlike the two previous methods, \Cref{algo:sxyz} does not form the second-moment matrices. It works directly on the data: it takes the SVD of $\mX$, projects $\mY$ onto the $r$ left singular vectors, truncates the projection to rank $k$ with a second SVD, and maps the result ($W$) back from the $r$ dimensional input space to the full input space. It returns $\mM$ only; the optional factorization into $(\mA, \mB)$ of \Cref{th:SXYZ} requires one more SVD and is not counted here. The total complexity is
$$ \min(n,p) \left[ n (C_{\text{svd}} p + q) + (C_{\text{svd}} \min(n, p, q) + p) q \right],$$
which for $p \approx q \approx k$ simplifies to $p \min(n, p)(C_{\text{svd}}n + n + C_{\text{svd}}\min(n,p) + p) \approx p \min(n, p)(10n + 10\min(n,p) + p)$. The dependence on $n$ is now through $\min(n,p)$. When $p \gg n$, the SVD of $\mX$ costs only $10n^2p$ and the total is dominated by the matrix products, $np^2$, which is much cheaper than the $p^3$ decompositions of the two other methods. When $n \gg p$ however, the SVD of $\mX$ costs $10np^2$, that is, about $10$ times the cost of computing the second-moment matrices.

\subsubsection*{Summary}

\Cref{tab:complexity} gathers the asymptotic complexities of the three methods for $p \approx q \approx k$ in the general case and in the two regimes $n \gg p$ and $p \gg n$. It shows that:
\begin{itemize}
    \item for $n \gg p$, \projected{} and \balanced{} are about $10$ times cheaper than \sxyz{}, because the second-moment matrices are cheaper to compute than the SVD of $\mX$;
    \item for $p \gg n$, \sxyz{} is the cheapest, because it avoids the $p^3$ decompositions of $\mS$;
    \item \projected{} and \balanced{} have the same asymptotic complexity, but \balanced{} relies on an SVD where \projected{} relies on an eigendecomposition.
\end{itemize}
These three predictions are confirmed by the fitting times measured on CPU and GPU in \Cref{sec:experiments}.

\begin{table}[ht]
\centering
\caption{Asymptotic time complexity of the three methods for $p \approx q \approx k$.}
\label{tab:complexity}
\begin{tabular}{l|ccc}
     & \textbf{Projected RRR} & \textbf{Balanced RRR} & \textbf{SXYZ RRR} \\
    \hline
    General & $np^2 + 20p^3$ & $np^2 + 20p^3$ & $10np^2 + 10p^3$ \\
    $n \gg p$ & $np^2$ & $np^2$ & $10np^2$ \\
    $p \gg n$ & $20p^3$ & $20p^3$ & $np^2$ \\
\end{tabular}
\end{table}

\section{Experiments}\label{sec:experiments}

\subsection{Protocol}

To investigate the practical performance of the algorithms presented in \Cref{sec:algorithms}, we conduct a series of experiments on synthetic datasets. We generated datasets of different sizes $n$ and considered different input dimensions $p$. We always choose the dimension of the output to be equal to the dimension of the input, i.e.\ $p=q$. For each setting $(n, p)$, we used three different methods to generate the input data and five different methods to generate the output data, leading to 15 different datasets for each choice of $(n, p)$. Then, we solve \cref{eq:defABstar} with the three introduced methods for two rank constraints, $k/p \in \{0.1, 0.9\}$. For each algorithm of \cref{sec:algorithms} and each dataset, we perform 10 runs with different random seeds; the seed also sets the generated data, so each run uses a new draw of $(\mX, \mY)$. Each method is tested on CPU and GPU with two float precisions (float32 and float64), except for $n=10^6$ in float64 on GPU, which was not run (the configuration $n=10^6$, $p=q=1000$ does not fit in the GPU memory). For more details on the experimental protocol, please refer to \Cref{detailed_experimental_protocol} and for more results refer to \Cref{experiments_extra}.


\subsection{Fitting time}

\begin{figure}
    \begin{center}
        \includegraphics[width=0.85\textwidth]{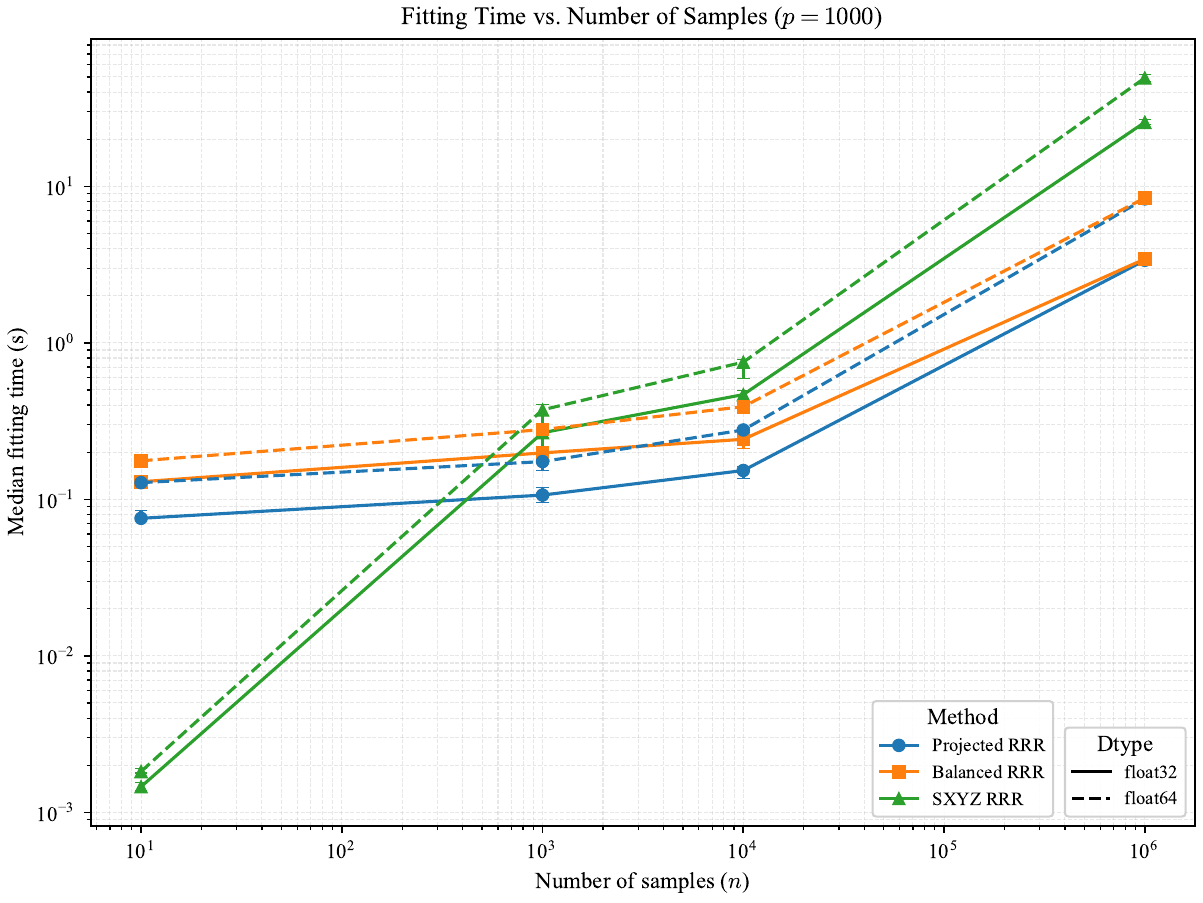}
    \end{center}
    \caption{Fitting time of the 3 algorithms on CPU for different $n$ and $p=q=1000$: median and 10th to 90th percentiles over data configurations, rank ratios and seeds.}
    \label{fig:fitting_time_cpu}
\end{figure}

\begin{figure}
    \begin{center}
        \includegraphics[width=0.85\textwidth]{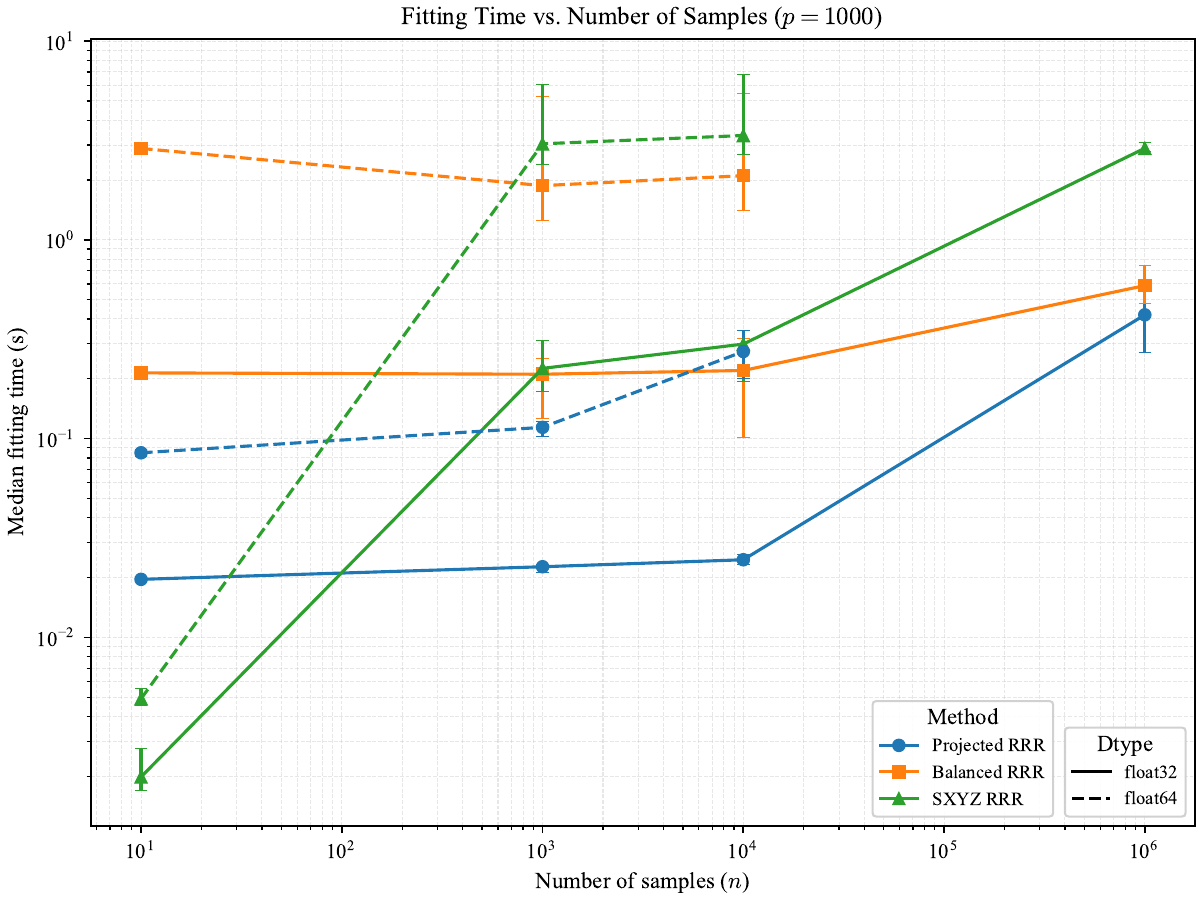}
    \end{center}
    \caption{Fitting time of the 3 algorithms on GPU for different $n$ and $p=q=1000$, aggregated as in \Cref{fig:fitting_time_cpu}. The float64 runs with $n=10^6$ were not performed.}
    \label{fig:fitting_time_gpu}
\end{figure}

\Cref{fig:fitting_time_cpu,fig:fitting_time_gpu} show the fitting times on CPU and GPU. The timings of \projected{} and \balanced{} include the computation of the factors $(\mA, \mB)$, while \sxyz{} only computes $\mM$.

On both devices, \sxyz{} is the fastest method for $n<p$ and \projected{} for $n \geq p$, in line with the complexity analysis of \Cref{sec:algorithms}. On GPU, \balanced{} is much slower than \projected{}, because its cost is dominated by an SVD rather than a eigendecomposition (which tend to be a few times faster). Float64 is slower than float32 on both devices, including for large $n$, since $\mS$ and $\mN$ are computed in the precision of the fit. Finally, the GPU pays off for large instances, but not always for small ones, in particular in float64. Note that the CPU and GPU runs were performed on different machines (\Cref{detailed_experimental_protocol}).

\subsection{Numerical performance}

We measure the numerical quality of the solutions with the mean squared error on the training data, $\frac{1}{nq} \norm{\mY - \mX \hat{\mM}}_F^2$. As the errors of different data configurations differ by orders of magnitude, each table reports the mean over the seeds for a single data configuration: the correlated-full-rank-noisy output with the isotropic input (\Cref{tab:evaluation_scores_cpu_isotropic_1000}, and \Cref{tab:evaluation_scores_cuda_isotropic_1000} in the appendix) or the rank-deficient anisotropic-low-rank input (\Cref{tab:evaluation_scores_cuda_anisotropic-low-rank_1000}, and \Cref{tab:evaluation_scores_cpu_anisotropic-low-rank_1000} in the appendix).

\begin{table}[htbp]
\centering
\caption{Training error on CPU, $\frac{1}{nq}\|\mathbf{Y} - \mathbf{X}\hat{\mathbf{M}}\|_F^2$, mean over the seeds (lower is better). Input: isotropic, output: correlated-full-rank-noisy, $p=1000$. Bold: lowest value of the row.}
\label{tab:evaluation_scores_cpu_isotropic_1000}
\begin{tabular}{cccccccc}
\toprule
$n$ & $k/p$ & \multicolumn{2}{c}{Projected} & \multicolumn{2}{c}{Balanced} & \multicolumn{2}{c}{SXYZ} \\
\cmidrule(lr){3-4}\cmidrule(lr){5-6}\cmidrule(lr){7-8}
 &  & float32 & float64 & float32 & float64 & float32 & float64 \\
\midrule
\multirow{2}{*}{1,000,000} & 0.1 & $\mathbf{5.54e3}$ & $\mathbf{5.54e3}$ & $\mathbf{5.54e3}$ & $\mathbf{5.54e3}$ & $\mathbf{5.54e3}$ & $\mathbf{5.54e3}$ \\
 & 0.9 & $\mathbf{4.98e-1}$ & $5.04e-1$ & $\mathbf{4.98e-1}$ & $5.04e-1$ & $\mathbf{4.98e-1}$ & $5.04e-1$ \\
\addlinespace
\multirow{2}{*}{10,000} & 0.1 & $5.43e3$ & $\mathbf{5.42e3}$ & $5.43e3$ & $\mathbf{5.42e3}$ & $5.43e3$ & $\mathbf{5.42e3}$ \\
 & 0.9 & $\mathbf{4.50e-1}$ & $4.56e-1$ & $\mathbf{4.50e-1}$ & $4.56e-1$ & $\mathbf{4.50e-1}$ & $4.56e-1$ \\
\addlinespace
\multirow{2}{*}{1,000} & 0.1 & $4.55e3$ & $\mathbf{4.54e3}$ & $4.55e3$ & $\mathbf{4.54e3}$ & $4.55e3$ & $\mathbf{4.54e3}$ \\
 & 0.9 & $5.58e-2$ & $\mathbf{3.25e-2}$ & $5.54e-2$ & $\mathbf{3.25e-2}$ & $3.26e-2$ & $\mathbf{3.25e-2}$ \\
\addlinespace
\multirow{2}{*}{10} & 0.1 & $1.06e-8$ & $\mathbf{2.04e-26}$ & $1.61e-8$ & $6.28e-26$ & $5.96e-9$ & $4.63e-26$ \\
 & 0.9 & $1.14e-8$ & $\mathbf{2.28e-26}$ & $1.61e-8$ & $6.28e-26$ & $5.96e-9$ & $4.63e-26$ \\
\bottomrule
\end{tabular}
\end{table}

\begin{table}[htbp]
\centering
\caption{Training error on GPU, $\frac{1}{nq}\|\mathbf{Y} - \mathbf{X}\hat{\mathbf{M}}\|_F^2$, mean over the seeds (lower is better). Input: anisotropic-low-rank, output: correlated-full-rank-noisy, $p=1000$. Bold: lowest value of the row. ---: not run.}
\label{tab:evaluation_scores_cuda_anisotropic-low-rank_1000}
\begin{tabular}{cccccccc}
\toprule
$n$ & $k/p$ & \multicolumn{2}{c}{Projected} & \multicolumn{2}{c}{Balanced} & \multicolumn{2}{c}{SXYZ} \\
\cmidrule(lr){3-4}\cmidrule(lr){5-6}\cmidrule(lr){7-8}
 &  & float32 & float64 & float32 & float64 & float32 & float64 \\
\midrule
\multirow{2}{*}{1,000,000} & 0.1 & $\mathbf{1.88e8}$ & --- & $\mathbf{1.88e8}$ & --- & $\mathbf{1.88e8}$ & --- \\
 & 0.9 & $\mathbf{4.87e-2}$ & --- & $3.28e1$ & --- & $5.60e0$ & --- \\
\addlinespace
\multirow{2}{*}{10,000} & 0.1 & $\mathbf{1.84e8}$ & $1.85e8$ & $\mathbf{1.84e8}$ & $1.85e8$ & $\mathbf{1.84e8}$ & $1.85e8$ \\
 & 0.9 & $4.16e-2$ & $\mathbf{3.80e-2}$ & $3.19e1$ & $\mathbf{3.80e-2}$ & $5.55e0$ & $\mathbf{3.80e-2}$ \\
\addlinespace
\multirow{2}{*}{1,000} & 0.1 & $\mathbf{1.56e8}$ & $1.57e8$ & $\mathbf{1.56e8}$ & $1.57e8$ & $\mathbf{1.56e8}$ & $1.57e8$ \\
 & 0.9 & $3.08e-2$ & $\mathbf{2.00e-2}$ & $3.02e1$ & $\mathbf{2.00e-2}$ & $3.90e0$ & $\mathbf{2.00e-2}$ \\
\addlinespace
\multirow{2}{*}{10} & 0.1 & $2.44e-4$ & $\mathbf{9.99e-22}$ & $7.03e-2$ & $1.92e-18$ & $1.56e-3$ & $5.81e-21$ \\
 & 0.9 & $2.77e-4$ & $\mathbf{1.58e-21}$ & $7.03e-2$ & $1.92e-18$ & $1.56e-3$ & $5.81e-21$ \\
\bottomrule
\end{tabular}
\end{table}

Most of the time, the three methods lead to similar errors. For $n=10<p$, the rank constraint is inactive and all the errors are rounding noise. The exceptions, mostly in float32 and especially on GPU, favor \projected{} or \sxyz{}:
\begin{itemize}
    \item for $n=p$ in float32 with the isotropic input, the errors of \projected{} and \balanced{} are about $70\%$ higher than the error of \sxyz{}, on both devices (\Cref{tab:evaluation_scores_cpu_isotropic_1000,tab:evaluation_scores_cuda_isotropic_1000});
    \item in float32 on GPU with the anisotropic-low-rank input, the errors of \balanced{} and \sxyz{} are two to three orders of magnitude higher than the error of \projected{} (\Cref{tab:evaluation_scores_cuda_anisotropic-low-rank_1000}).
\end{itemize}

\subsection{Practical recommendations}

\begin{itemize}
    \item \textbf{Algorithm.} Use \sxyz{} for $n<p$ and \projected{} for $n \geq p$. For large $n$, if memory is a concern, \projected{} can be run in batches.
    \item \textbf{Precision.} float32 is usually enough, but when $\mX$ may be ill conditioned or rank deficient, use float64, with which the three methods almost always reach the same error. It costs up to $2.5\times$ more time on CPU and $14\times$ on GPU in our runs.
    \item \textbf{Device.} Use a GPU for large instances: in float32 with $n=10^6$, it was $6$ to $9\times$ faster than our CPU. For small instances, the CPU can be as fast or faster, in particular in float64.
\end{itemize}

\bibliographystyle{plainnat}
\bibliography{bibli}

\appendix
\newpage

\section{Proofs}
\begin{theorem}
    \thbalanced
\end{theorem}
\begin{proof}\label{proof:defOLSequi}
We solve the minimization problem defined in \cref{eq:defABstar}, which we recall here:
\begin{align*}
    \mA^*, \mB^* \in \argmin_{\mA, \mB}\norm{\mX\mA\mB - \mY}^2 \;\; s.t. \;\; \mA \in \R^{p,k}, \mB \in \R^{k, q}
\end{align*}

First, we show that:
\begin{equation}
    \argmin_{\mA, \mB}\norm{\mX\mA\mB - \mY} = \argmin_{\mA, \mB}\norm{\Shalf\mA\mB - \Smhalf\mN}
\end{equation}
where $\Shalf$ and $\Smhalf$ are the square root and the pseudo-inverse of the square root of $\mS = \frac{1}{n}\trans{\mX}\mX$, and $\mN = \frac{1}{n} \trans{\mX}\mY$. Indeed, we have:

\begin{equation}
    \frac{1}{n} \norm{\mX\mA\mB - \mY}^2 - \frac{1}{n}\norm{\mY}^2 = \norm{\Shalf\mA\mB - \Smhalf\mN}^2 - \norm{\Smhalf\mN}^2\label{eq:step1}
\end{equation}
\begin{myproof}
\underline{Proof of \cref{eq:step1}}:
Let $\mA \in \R^{p,k}, \mB \in \R^{k, q}$:
\begin{align}
    \frac{1}{n}\norm{\mX\mA\mB - \mY}^2 - \frac{1}{n}\norm{\mY}^2
    &= \frac{1}{n}\norm{\mX\mA\mB}^2 -\frac{2}{n}\scalarp{\mX\mA\mB ,\; \mY} + \frac{1}{n}\norm{\mY}^2 - \frac{1}{n}\norm{\mY}^2\\
    &=\trace{\trans{\mB}\trans{\mA}\frac{\trans{\mX}\mX}{n}\mA\mB} -2\trace{\trans{\mB}\trans{\mA}\frac{\trans{\mX}\mY}{n}} \\
    &= \trace{\trans{\mB} \trans{\mA} \mS \mA\mB} -2\trace{\trans{\mB}\trans{\mA} \mN}
\intertext{
As $\Shalf\Smhalf\mN=\mN$ (see \cref{lem:SSN=N}) and $\Shalf=\trans{\Shalf}$, it follows that:}
    \frac{1}{n}\norm{\mX\mA\mB - \mY}^2 - \frac{1}{n}\norm{\mY}^2&=\trace{\trans{\mB}\trans{\mA}\trans{\Shalf}\Shalf\mA\mB} -2\trace{\trans{\mB}\trans{\mA}\trans{\Shalf}\Smhalf\mN}\\
     &= \norm{\Shalf\mA\mB}^2 -2\scalarp{\Shalf\mA\mB, \Smhalf\mN} \\
     &= \norm{\Shalf\mA\mB - \Smhalf\mN}^2 -\norm{\Smhalf\mN}^2
\end{align}
\end{myproof}

\underline{Now that we have proved \cref{eq:step1}, we minimize its right-hand side over $(\mA, \mB)$}:

Let
\begin{align}
    \mathcal{U} := \argmin_{\mA, \mB}\norm{\Shalf\mA\mB - \Smhalf\mN}^2\;\; \suchthat \;\; \mA \in \R^{p,k}, \mB \in \R^{k, q}
\end{align}
Then, using \cref{lem:AInSA}, the set of solutions
\begin{align}
    \mathcal{V} := \left\{\left(\Smhalf\mA_v, \mB_v\right)\;| \; \left(\mA_v, \mB_v\right) \in\argmin_{\mA, \mB}\norm{\mA\mB - \Smhalf\mN}\;\;\suchthat\;\; \mA \in \R^{p,k}, \mB \in \R^{k, q}\right\}
\end{align}
is included in $\mathcal{U}$, and we look for a solution in $\mathcal{V}$. A pair $\left(\mA_v, \mB_v\right)$ such that $\left(\Smhalf\mA_v, \mB_v\right)$ is in $\mathcal{V}$ is given by \citet{ApproxMatrix}: with $\trans{\mU}\Lambda\mV$ the reduced SVD of $\Smhalf\mN$, a solution is
\begin{align*}
     \mA_v &= \trans{\mU_{:k}}\sqrt{\Lambda_{:k, :k}}\\
        \mB_v &= \sqrt{\Lambda_{:k, :k}} \mV_{:k}
\end{align*}
Hence, $(\mA^*, \mB^*) = (\Smhalf\mA_v, \mB_v)$ is a solution of \cref{eq:defABstar}, and the induced solution of \cref{eq:Mstar} is
\begin{align*}
    \mA^*\mB^* = \Smhalf\trans{\mU_{:k}}\Lambda_{:k, :k}\mV_{:k} = \Smhalf\left[\Smhalf\mN\right]_k = \mM^*_{\sqrt{\;}}\;.
\end{align*}
\end{proof}
\vspace{1cm}
\begin{theorem}
    \thSXYZ
\end{theorem}
\begin{proof}\label{proof:SXYZ}
We solve the minimization problem defined in \cref{eq:Mstar}, which we recall here:
\begin{align}
    \mM^* \in \argmin_{\mM\in\mathbb{R}^{p, q}}\norm{\mX\mM - \mY} \;\; s.t. \;\; \rank(\mM) \leq k \label{eq:Mstarxx}
\end{align} 
Let $\mX = \trans{\tilde{\mE}}\tilde{\mD}\tilde{\mF}$ be the \textbf{full} SVD of $\mX\in\mathbb{R}^{n, p}$. We have $\tilde{\mE}\in\mathbb{R}^{n, n}$, $\tilde{\mD}\in\mathbb{R}^{n, p}$ and  $\tilde{\mF}\in\mathbb{R}^{p, p}$ where the matrices $\tilde{\mE}$ and $\tilde{\mF}$ are orthogonal. 
Because $\tilde{\mE}$ is an orthogonal matrix:
\begin{align}
    \norm{\mX\mM - \mY}
    &= \norm{\tilde{\mE}(\mX\mM - \mY)}\\
    &=\norm{\tilde{\mD}\tilde{\mF}\mM - \tilde{\mE}\mY}
\end{align}
As $\tilde{\mF}$ is invertible, we can perform the change of variable $\tilde{\mM} = \tilde{\mF}\mM$. Therefore, solving
\begin{equation}\label{eq:Mstarxxbis}
    \tilde{\mM}^* \in \argmin_{\tilde{\mM}\in\mathbb{R}^{p, q}}\norm{\tilde{\mD}\tilde{\mM} - \tilde{\mE}\mY} \;\; s.t. \;\; \rank(\tilde{\mM}) \leq k
\end{equation}
gives a solution of \cref{eq:Mstarxx} with $\mM^* = {\tilde{\mF}}^\top\tilde{\mM}^*$.
Let $r=\rank(\mX)$. We write $\tilde{\mD} = \begin{pmatrix}
    \mD &\mathbf{0}\\
    \mathbf{0} & \mathbf{0}
\end{pmatrix}$, $\tilde{\mE}= \begin{pmatrix}
    \mE\\
    \mE^{\perp}
\end{pmatrix}$ and $\tilde{\mF}= \begin{pmatrix}
    \mF\\
    \mF^{\perp}
\end{pmatrix}$, where $\mD \in\mathbb{R}^{r, r}$, $\mE\in\mathbb{R}^{r, n}$ and $\mF\in\mathbb{R}^{r, p}$ correspond to the reduced SVD of $\mX$. For any matrix $\tilde{\mM}\in\mathbb{R}^{p, q}$, we have
\begin{align}
    \norm{\tilde{\mD}\tilde{\mM} - \tilde{\mE}\mY}^2 &= \norm{\begin{pmatrix}
    \mD &\mathbf{0}\\
    \mathbf{0} & \mathbf{0}
\end{pmatrix}\tilde{\mM} - \begin{pmatrix}
    \mE\\
    \mE^{\perp}
\end{pmatrix}\mY}^2\\
&= \norm{\mD\tilde{\mM}_{:r, :} - \mE\mY}^2 + \norm{\mE^{\perp}\mY}^2
\end{align}
The matrix $\mD$ being invertible, we perform the change of variable $\mW = \mD\tilde{\mM}_{:r, :}\in\mathbb{R}^{r, q}$ and solve the minimization problem:
\begin{equation}
    \mW^* \in \argmin_{\mW\in\mathbb{R}^{r, q}}\norm{\mW - \mE\mY} \;\; \text{s.t.} \;\; \rank(\mW) \leq k
\end{equation}
By \citet{ApproxMatrix}, a solution of this problem is $\mW^* = \left[\mE\mY\right]_k$. This problem gives a lower bound for \cref{eq:Mstarxxbis}: the objective of \cref{eq:Mstarxxbis} does not depend on the last $p-r$ rows of $\tilde{\mM}$, and for any feasible $\tilde{\mM}$, the matrix $\mW = \mD\tilde{\mM}_{:r, :}$ is feasible, since $\rank(\mW) = \rank(\tilde{\mM}_{:r, :}) \leq \rank(\tilde{\mM}) \leq k$. The matrix $\tilde{\mM}^*$ with $\tilde{\mM}^*_{:r, :} = \mD^{-1}\mW^*$ and zeros in its last $p-r$ rows satisfies $\rank(\tilde{\mM}^*) = \rank(\mW^*) \leq k$ and attains this bound, so it is a solution of \cref{eq:Mstarxxbis}. The associated $\mM^*$ is
\begin{align}
    \mM^* &= \trans{\tilde{\mF}}\tilde{\mM}^*\\
    &=\begin{pmatrix}
        \trans{\mF} & \trans{{\mF^{\perp}}}
    \end{pmatrix} \begin{pmatrix}
        \mD^{-1}[\mE\mY]_k\\
        \mathbf{0}
    \end{pmatrix}\\
    &=\trans{\mF}\mD^{-1}[\mE\mY]_k \;.
\end{align}
\end{proof}
\vspace{1cm}

\begin{proposition}
    \propMsqrtequaldMperp
\end{proposition}
\begin{proof}\label{prop:Msqrt=Mperpproof}
    By the common choice of \cref{rem:2}, the reduced SVD of $\Smhalf\mN$ and the one of $\trans{\mN}\mS^+\mN$ share the same matrix of right singular vectors, so that we can write $\Smhalf\mN = \trans{\mU}\Lambda\mV$ and $\trans{\mN}\mS^+\mN = \trans{\mV}\Lambda^2\mV$ with $r=\rank(\Smhalf\mN)$, $\mU\in\mathbb{R}^{r, p}$, $\mV\in\mathbb{R}^{r, q}$ and $\Lambda\in\mathbb{R}^{r, r}$.
    It follows that:
    \begin{align}
        \mM^*_{\perp} &= \mS^+\mN\trans{\mV_{:k}}\mV_{:k} \\
            &= \Smhalf\Smhalf\mN\trans{\mV_{:k}}\mV_{:k}\\
            &= \Smhalf\trans{\mU}\Lambda \mV\trans{\mV_{:k}}\mV_{:k} \\
            &= \Smhalf\trans{\mU}\Lambda \begin{pmatrix}
                \mV_{:k} \\ \mV_{k+1:}
            \end{pmatrix}\trans{\mV_{:k}}\mV_{:k}
    \end{align}
    Because the rows of $\mV$ are orthonormal, we have $\mV_{:k}\trans{\mV_{:k}} = \mathbf{I}_k \in\mathbb{R}^{k, k}$ with $\mI_k$ the identity matrix, and $\mV_{k+1:}\trans{\mV_{:k}} = \mathbf{0}_{r-k, k}\in\mathbb{R}^{r-k, k}$, with $\mathbf{0}_{k, k'}$ the $k \times k'$ matrix full of zeros. It follows that:
    \begin{align}
        \mM^*_{\perp} &= \Smhalf\trans{\mU}\Lambda \begin{pmatrix}
                \mathbf{I}_k \\ \mathbf{0}_{r-k, k}
            \end{pmatrix}\mV_{:k}\\
            &= \Smhalf\trans{\mU}\begin{pmatrix}
        \Lambda_{:k, :k} & \mathbf{0}_{k, r-k} \\
        \mathbf{0}_{r-k, k} & \Lambda_{k+1:, k+1:}
    \end{pmatrix} \begin{pmatrix}
                \mathbf{I}_k \\ \mathbf{0}_{r-k, k}
            \end{pmatrix}\mV_{:k}\\
    &=\Smhalf \begin{pmatrix}
        \mU_{:k} \\
        \mU_{k+1:}
    \end{pmatrix}^{\top}
    \begin{pmatrix}
        \Lambda_{:k, :k} \\
        \mathbf{0}_{r-k, k} 
    \end{pmatrix} \mV_{:k}\\
    &=\Smhalf \trans{\mU_{:k}} \Lambda_{:k, :k} \mV_{:k}\\
    &= \Smhalf\left[\Smhalf\mN\right]_k \\
        &= \mM^*_{\sqrt{\;}}
    \end{align}
\end{proof}
\begin{proposition}
    \propMperpequalMcirc
\end{proposition}
\begin{proof}\label{prop:circ=Mperpproof}
    For $\mX\in\R^{n, p}$ and $\mY\in\R^{n, q}$, we recall the following notations.\\
    \begin{minipage}{0.48\textwidth}
        \begin{align*}
            \mS &= \frac{1}{n}\trans{\mX}\mX\in\R^{p, p}\\
            \mN &= \frac{1}{n}\trans{\mX}\mY\in\R^{p, q}\\
            \mX &= \trans{\mE}\mD\mF &&\text{reduced SVD}\\
            \trans{\mN}\mS^+\mN &= \trans{\mV}\Lambda^2\mV  &&\text{reduced SVD}
        \end{align*}
    \end{minipage}
    \hfill
    \begin{minipage}{0.48\textwidth}
        \begin{align*}
            \mM^*_{\perp} &= \mS^+\mN\trans{\mV_{:k}}\mV_{:k}    \\
            \mM^*_{\circ} &= \trans{\mF}\mD^{-1}[\mE\mY]_k
        \end{align*}
    \end{minipage}
    We now use the reduced SVD of $\mX$ to express $\mM^*_{\perp}$. The factors $\frac{1}{n}$ cancel in $\mS^+\mN = \big(\trans{\mX}\mX\big)^{+}\trans{\mX}\mY$, so that:
    \begin{align}
        \mM^*_{\perp} &= \big(\trans{\mX}\mX\big)^{+}\trans{\mX}\mY\trans{\mV_{:k}}\mV_{:k}\\
        &=\big(\trans{\mF}\mD^{2}\mF\big)^{+}\trans{\mF}\mD\mE\mY\trans{\mV_{:k}}\mV_{:k}\\
        &=\trans{\mF}\mD^{-1}\mE\mY\trans{\mV_{:k}}\mV_{:k}\label{eq:MperAlmostMcirc}
    \end{align}
    To conclude, we show that $\mE\mY\trans{\mV_{:k}}\mV_{:k} = [\mE\mY]_k$. Using the SVD of $\mX$ to express $\trans{\mN}\mS^+\mN$, we have:
    \begin{align}
        \trans{\mN}\mS^+\mN &= \frac{1}{n}\trans{\mY}\mX\big(\trans{\mX}\mX\big)^{+}\trans{\mX}\mY \\
        &= \frac{1}{n}\trans{\mY}\trans{\mE}\mD\mF\big(\trans{\mF}\mD\mE\trans{\mE}\mD\mF\big)^{+}\trans{\mF}\mD\mE\mY\\
        &= \frac{1}{n}\trans{\mY}\trans{\mE}\mE\mY
    \end{align}
    Hence $\trans{\mN}\mS^+\mN = \frac{1}{n}\trans{(\mE\mY)}(\mE\mY)$, so the singular values of $\mE\mY$ are $\Lambda' := \sqrt{n}\Lambda$ and, by the common choice of \cref{rem:2}, its right singular vectors are the matrix $\mV$ above. In other words, the reduced SVD of $\mE\mY$ is
    \begin{align}
        \mE\mY = \trans{\mW}\Lambda'\mV\;,
    \end{align}
    As a consequence, for $w_j$ and $v_j$ the $j^{th}$ left and right singular vectors of $\mE\mY$, i.e. the $j^{th}$ rows of $\mW$ and $\mV$, we have
    \begin{align}
        \mE\mY\trans{\mV_{:k}}\mV_{:k} &= \sum_{j}w_j\trans{v_j}\Lambda'_{j, j}\sum_{f=1}^kv_f\trans{v_f}\\
        &= \sum_{j=1}^kw_j\trans{v_j}\Lambda'_{j, j}\\
        &= [\mE\mY]_k\;.
    \end{align}
    Continuing from \cref{eq:MperAlmostMcirc}, it follows that
    \begin{align}
        \mM^*_{\perp} &= \trans{\mF}\mD^{-1}[\mE\mY]_k\\
        &=\mM^*_{\circ} \;.
    \end{align}
\end{proof}

\begin{lemma}\label{lem:SSN=N}
Let $\mX \in \mathbb{R}^{n, p}$ and $\mY \in \mathbb{R}^{n, q}$, and define $\mS := \trans{\mX}\mX$ and $\mN :=\trans{\mX}\mY$. Then
\begin{align}
    \Shalf\Smhalf\mN &= \mN
\end{align}
\end{lemma}
\begin{proof}
For $\mX\in\mathbb{R}^{n, p}$, let $\trans{\mE}\mD\mF$ be its reduced SVD with $r = \rank(\mX)$, $\mE \in \mathbb{R}^{r, n}$, $\mF \in \mathbb{R}^{r, p}$ and $\mD \in \mathbb{R}^{r, r}$ with the elements of its  diagonal being strictly positive. We have $\mS = \trans{\mF}\mD^2\mF$ and 
    \begin{align*}
        \Shalf\Smhalf\mN &= \trans{\mF}\sqrt{\mD^2}\mF\trans{\mF}\sqrt{\mD^2}^{+}\mF\;\mN
    \end{align*}
We have $\mF\trans{\mF} = \mI \in \mathbb{R}^{r, r}$, and $\mD^2$ is invertible, so $\sqrt{\mD^2}\sqrt{\mD^2}^{+} = \mI$. It follows that 
\begin{align*}
     \Shalf\Smhalf\mN &= \trans{\mF}\mF\mN\\
     &=\trans{\mF}\mF\trans{\mX}\mY\\
     &=\trans{\mF}\mF\trans{\mF}\mD\mE\mY \\ 
     &=\trans{\mF}\mD\mE\mY \\ 
     &= \trans{\mX}\mY\\
     &=\mN
\end{align*}
\end{proof}

\begin{lemma}\label{lem:AInSA}
Let $\mS\in\mathbb{R}^{p, p}$ be a symmetric positive semi-definite matrix, $\mN\in\mathbb{R}^{p, q}$ and $k\in\mathbb{N}^*$. Let
    \begin{align}
    \mathcal{U} &:=\argmin_{\mA, \mB}\norm{\Shalf\mA\mB - \Smhalf\mN}\;\;\suchthat\;\; \mA \in \R^{p,k}, \mB \in \R^{k, q} \label{eq:minAS1demi}\\
    \mathcal{V} &:= \left\{\left(\Smhalf\mA, \mB\right)\;| \; \left(\mA, \mB\right) \in\argmin_{\mA, \mB}\norm{\mA\mB - \Smhalf\mN}\;\;\suchthat\;\; \mA \in \R^{p,k}, \mB \in \R^{k, q}\right\}\label{eq:minA}
\end{align}
Then $\mathcal{V} \subset \mathcal{U}$.
\end{lemma}
\begin{proof}
Let $\left(\Smhalf\mA, \mB\right)\in \mathcal{V}$, and let $\mA_s \in \R^{p,k}$, $\mB_s \in \R^{k, q}$ be arbitrary. We show that:
\begin{align}
\norm{\Shalf\Smhalf\mA\mB - \Smhalf\mN} \leq \norm{\Shalf\mA_s\mB_s - \Smhalf\mN}
\end{align}
Because $\mathbb{R}^p = \Ima(\Smhalf)\oplus_{\perp}\ker(\Shalf)$ and $\ker(\Shalf) = \ker(\Smhalf)$, let $\mV := \begin{pmatrix} \vv_1, ...\vv_k\end{pmatrix}$ and $\mU := \begin{pmatrix}\vu_1, ..., \vu_k \end{pmatrix}$, with $\vu_i\in\ker(\Smhalf)$ for all $i=1, ..., k$, be such that
\begin{align}
    \mA = \Smhalf\mV + \mU
\end{align}
It follows that 
\begin{align*}
    \norm{\mA\mB - \Smhalf\mN}^2 &= \norm{ \Smhalf\mV\mB + \mU\mB- \Smhalf\mN}^2 \\
    &= \norm{ \Smhalf\mV\mB - \Smhalf\mN}^2 + 2 \scalarp{\mU\mB, \Smhalf\left(\mV\mB - \mN\right)} + \norm{\mU\mB}^2
\end{align*}
Because $\vu_i\in\ker(\Smhalf)$ and $\Smhalf$ is symmetric, we have $\trans{\vu_i}\Smhalf = \Smhalf\vu_i = 0$ and $\mU^\top\Smhalf = \mathbf{0}$. It follows that:
\begin{align}
    \norm{\mA\mB - \Smhalf\mN}^2 &= \norm{ \Smhalf\mV\mB - \Smhalf\mN}^2 + \norm{\mU\mB}^2 \label{eq:defAB}
\end{align}
The pair $(\Smhalf\mV, \mB)$ is also feasible for the problem defining $\mathcal{V}$, and by \cref{eq:defAB}, its objective is smaller than the one of $(\mA, \mB)$ by $\norm{\mU\mB}^2$. Since $(\mA, \mB)$ is a minimizer of this problem, $\mU\mB =\mathbf{0}$, i.e.\ $\mA\mB = \Smhalf\mV\mB$. Moreover, as $\Shalf\Smhalf\Smhalf = \Smhalf$ and $\Shalf\Smhalf\mU = \mathbf{0}$, we have:
\begin{align}
    \Shalf\Smhalf\mA\mB = \Smhalf\mV\mB = \mA\mB
\end{align}
Then
\begin{align}
    \norm{\Shalf\Smhalf\mA\mB - \Smhalf\mN} &= \norm{\mA\mB - \Smhalf\mN}\label{eq:last}
\end{align}
Because $(\Smhalf\mA, \mB)\in \mathcal{V}$, the pair $(\mA, \mB)$ is a minimizer of the problem defining $\mathcal{V}$, for which the pair $(\Shalf\mA_s, \mB_s)$ is feasible. Hence:
\begin{align}
     \norm{\mA\mB - \Smhalf\mN}&\leq\norm{\Shalf\mA_s\mB_s - \Smhalf\mN}\label{eq:lastlast}
\end{align}
Combining \Cref{eq:last,eq:lastlast}, it follows that
\begin{align}
    \norm{\Shalf\Smhalf\mA\mB - \Smhalf\mN} \leq \norm{\Shalf\mA_s\mB_s - \Smhalf\mN}
\end{align}
Since $(\mA_s, \mB_s)$ is arbitrary, $(\Smhalf\mA, \mB)\in\mathcal{U}$.
\end{proof}

\section{Detailed experimental protocol}

\label{detailed_experimental_protocol}

We benchmark the three algorithms (\sxyz{}, \projected{}, and \balanced{}) on synthetic datasets with controlled properties. All experiments are repeated over 10 random seeds.

\subsection{Data generation}

\paragraph{Input data $\mX \in \R^{n \times p}$.} We consider three generation schemes:
\begin{itemize}
    \item \textbf{Isotropic:} $\mX$ has i.i.d.\ entries drawn from $\N(0, \sigma^2)$ with $\sigma = 0.1$.
    \item \textbf{Anisotropic full-rank:} $\mX = \mS\mW$ where $\mS \in \R^{n \times \kappa}$, $\mW \in \R^{\kappa \times p}$ with i.i.d.\ standard normal entries and $\kappa = 2p \geq p$.
    \item \textbf{Anisotropic low-rank:} Same construction with $\kappa = \lfloor p/2 \rfloor < p$, so that $\rank(\mX) \leq \kappa$.
\end{itemize}

\paragraph{Output data $\mY \in \R^{n \times q}$.} For each input $\mX$, we generate $\mY$ according to five schemes:
\begin{itemize}
    \item \textbf{Independent:} $\mY$ has i.i.d.\ entries drawn from $\N(0, \lambda^2)$ with $\lambda = 0.2$, independently of $\mX$.
    \item \textbf{Correlated, low-rank:} $\mY = \mX\mA\mB$ where $\mA \in \R^{p \times r}$, $\mB \in \R^{r \times q}$ with i.i.d.\ standard normal entries and $r = \lfloor q/2 \rfloor$.
    \item \textbf{Correlated, full-rank:} Same as above with $r = q$.
    \item \textbf{Correlated with noise, low-rank:} $\mY = \mX\mA\mB + \vepsilon$ with $\vepsilon \sim \N(0, \lambda^2 \mI)$, $\lambda = 0.2$, and $r = \lfloor q/2 \rfloor$.
    \item \textbf{Correlated with noise, full-rank:} Same as above with $r = q$.
\end{itemize}

This yields $3 \times 5 = 15$ different data generation configurations for each choice of dimensions.

\subsection{Hyperparameters}

The experiments sweep over the following parameter grid:
\begin{itemize}
    \item Number of samples: $n \in \{10,\; 10^3,\; 10^4,\; 10^6\}$
    \item Feature dimensions: $(p, q) \in \{(10, 10),\; (100, 100),\; (1000, 1000)\}$
    \item Desired rank: $k \in \{0.1p,\; 0.9p\}$, giving $k = \max(\lfloor 0.1\,p \rfloor, 1)$ or $k = \lfloor 0.9\,p \rfloor$
\end{itemize}
All methods use the default spectral cut-off.
Each configuration is run on both CPU and CUDA, and with both \texttt{float32} and \texttt{float64} precision, except for $n=10^6$ in \texttt{float64} on CUDA, which was not run (the configuration $n=10^6$, $p=q=1000$ does not fit in the GPU memory).

\subsection{Metrics}

For each run, we record:
\begin{itemize}
    \item \textbf{Fitting time:} wall-clock time of the call that fits the model. For \projected{} and \balanced{}, it includes the computation of $\mS$ and $\mN$ (in a single batch) and of $\mA$, $\mB$ and $\mM$. For \sxyz{}, it only includes the computation of $\mM$: the factors $(\mA, \mB)$ are not computed.
    \item \textbf{Evaluation score:} the mean squared error $\frac{1}{nq}\norm{\mY - \mX\hat{\mM}}_F^2$ where $\hat{\mM}$ is the estimated coefficient matrix. It is computed on the training data, i.e.\ the data $(\mX, \mY)$ used for the fit.
\end{itemize}

\paragraph{Aggregation.} Each point of the fitting time figures aggregates all the runs with the same device, precision, method, $p$ and $n$, i.e.\ 15 data configurations $\times$ 2 rank ratios $\times$ 10 seeds $= 300$ runs. We report the median of these runs, and the error bars span their 10th to 90th percentiles. The error bars therefore mix different data configurations and rank ratios, and do not measure the timing noise of a single configuration. Each cell of the evaluation score tables is the mean over the 10 seeds of a single data configuration, rank ratio, $n$, device, precision and method.

\subsection{Hardware and software}

Experiments with GPU were run on:
\begin{itemize}
    \item CPU: Intel Xeon Silver 4214R @ 2.40GHz
    \item GPU: Nvidia A40 48G
\end{itemize}

Experiments without GPU were run on:
\begin{itemize}
    \item CPU: Intel Xeon E5-2650 v4 @ 2.20GHz
\end{itemize}

The CPU and GPU timings therefore come from different machines. Each benchmark ran as a Slurm job with 2 CPU cores allocated (\texttt{-{}-cpus-per-task=2}) and exclusive access to its node (\texttt{-{}-exclusive}), so that no other job ran on the same machine during the timings. All the runs used PyTorch 2.13.0 (the PyPI wheel pinned in the \texttt{uv.lock} file of the code), whose CPU linear algebra relies on Intel oneMKL 2024.2. On GPU, it uses CUDA 13.0, and the default heuristic of PyTorch picks between cuSOLVER and MAGMA 2.6.1 for each decomposition.

\section{More experimental results}\label{experiments_extra}

\begin{table}[htbp]
\centering
\caption{Training error on GPU, $\frac{1}{nq}\|\mathbf{Y} - \mathbf{X}\hat{\mathbf{M}}\|_F^2$, mean over the seeds (lower is better). Input: isotropic, output: correlated-full-rank-noisy, $p=1000$. Bold: lowest value of the row. ---: not run.}
\label{tab:evaluation_scores_cuda_isotropic_1000}
\begin{tabular}{cccccccc}
\toprule
$n$ & $k/p$ & \multicolumn{2}{c}{Projected} & \multicolumn{2}{c}{Balanced} & \multicolumn{2}{c}{SXYZ} \\
\cmidrule(lr){3-4}\cmidrule(lr){5-6}\cmidrule(lr){7-8}
 &  & float32 & float64 & float32 & float64 & float32 & float64 \\
\midrule
\multirow{2}{*}{1,000,000} & 0.1 & $\mathbf{5.53e3}$ & --- & $\mathbf{5.53e3}$ & --- & $\mathbf{5.53e3}$ & --- \\
 & 0.9 & $\mathbf{5.04e-1}$ & --- & $5.05e-1$ & --- & $5.05e-1$ & --- \\
\addlinespace
\multirow{2}{*}{10,000} & 0.1 & $\mathbf{5.41e3}$ & $5.42e3$ & $\mathbf{5.41e3}$ & $5.42e3$ & $\mathbf{5.41e3}$ & $5.42e3$ \\
 & 0.9 & $4.56e-1$ & $\mathbf{4.53e-1}$ & $4.57e-1$ & $\mathbf{4.53e-1}$ & $4.57e-1$ & $\mathbf{4.53e-1}$ \\
\addlinespace
\multirow{2}{*}{1,000} & 0.1 & $\mathbf{4.54e3}$ & $4.55e3$ & $\mathbf{4.54e3}$ & $4.55e3$ & $\mathbf{4.54e3}$ & $4.55e3$ \\
 & 0.9 & $5.55e-2$ & $\mathbf{3.31e-2}$ & $5.54e-2$ & $\mathbf{3.31e-2}$ & $3.33e-2$ & $\mathbf{3.31e-2}$ \\
\addlinespace
\multirow{2}{*}{10} & 0.1 & $5.36e-9$ & $\mathbf{2.19e-26}$ & $1.69e-6$ & $3.79e-23$ & $3.70e-8$ & $1.22e-25$ \\
 & 0.9 & $6.16e-9$ & $\mathbf{3.36e-26}$ & $1.69e-6$ & $3.79e-23$ & $3.70e-8$ & $1.22e-25$ \\
\bottomrule
\end{tabular}
\end{table}

\begin{table}[htbp]
\centering
\caption{Training error on CPU, $\frac{1}{nq}\|\mathbf{Y} - \mathbf{X}\hat{\mathbf{M}}\|_F^2$, mean over the seeds (lower is better). Input: anisotropic-low-rank, output: correlated-full-rank-noisy, $p=1000$. Bold: lowest value of the row.}
\label{tab:evaluation_scores_cpu_anisotropic-low-rank_1000}
\begin{tabular}{cccccccc}
\toprule
$n$ & $k/p$ & \multicolumn{2}{c}{Projected} & \multicolumn{2}{c}{Balanced} & \multicolumn{2}{c}{SXYZ} \\
\cmidrule(lr){3-4}\cmidrule(lr){5-6}\cmidrule(lr){7-8}
 &  & float32 & float64 & float32 & float64 & float32 & float64 \\
\midrule
\multirow{2}{*}{1,000,000} & 0.1 & $1.89e8$ & $\mathbf{1.88e8}$ & $1.89e8$ & $\mathbf{1.88e8}$ & $1.89e8$ & $\mathbf{1.88e8}$ \\
 & 0.9 & $4.35e-2$ & $\mathbf{4.00e-2}$ & $4.55e-2$ & $\mathbf{4.00e-2}$ & $4.32e-2$ & $\mathbf{4.00e-2}$ \\
\addlinespace
\multirow{2}{*}{10,000} & 0.1 & $\mathbf{1.85e8}$ & $\mathbf{1.85e8}$ & $\mathbf{1.85e8}$ & $\mathbf{1.85e8}$ & $\mathbf{1.85e8}$ & $\mathbf{1.85e8}$ \\
 & 0.9 & $4.04e-2$ & $\mathbf{3.80e-2}$ & $4.23e-2$ & $\mathbf{3.80e-2}$ & $4.05e-2$ & $\mathbf{3.80e-2}$ \\
\addlinespace
\multirow{2}{*}{1,000} & 0.1 & $1.57e8$ & $\mathbf{1.56e8}$ & $1.57e8$ & $\mathbf{1.56e8}$ & $1.57e8$ & $\mathbf{1.56e8}$ \\
 & 0.9 & $2.48e-2$ & $\mathbf{2.00e-2}$ & $2.53e-2$ & $\mathbf{2.00e-2}$ & $2.19e-2$ & $\mathbf{2.00e-2}$ \\
\addlinespace
\multirow{2}{*}{10} & 0.1 & $5.03e-4$ & $\mathbf{9.64e-22}$ & $8.96e-4$ & $3.67e-21$ & $2.69e-4$ & $3.54e-21$ \\
 & 0.9 & $5.48e-4$ & $\mathbf{1.09e-21}$ & $8.96e-4$ & $3.67e-21$ & $2.69e-4$ & $3.54e-21$ \\
\bottomrule
\end{tabular}
\end{table}

\begin{figure}
    \begin{center}
        \includegraphics[width=0.9\textwidth]{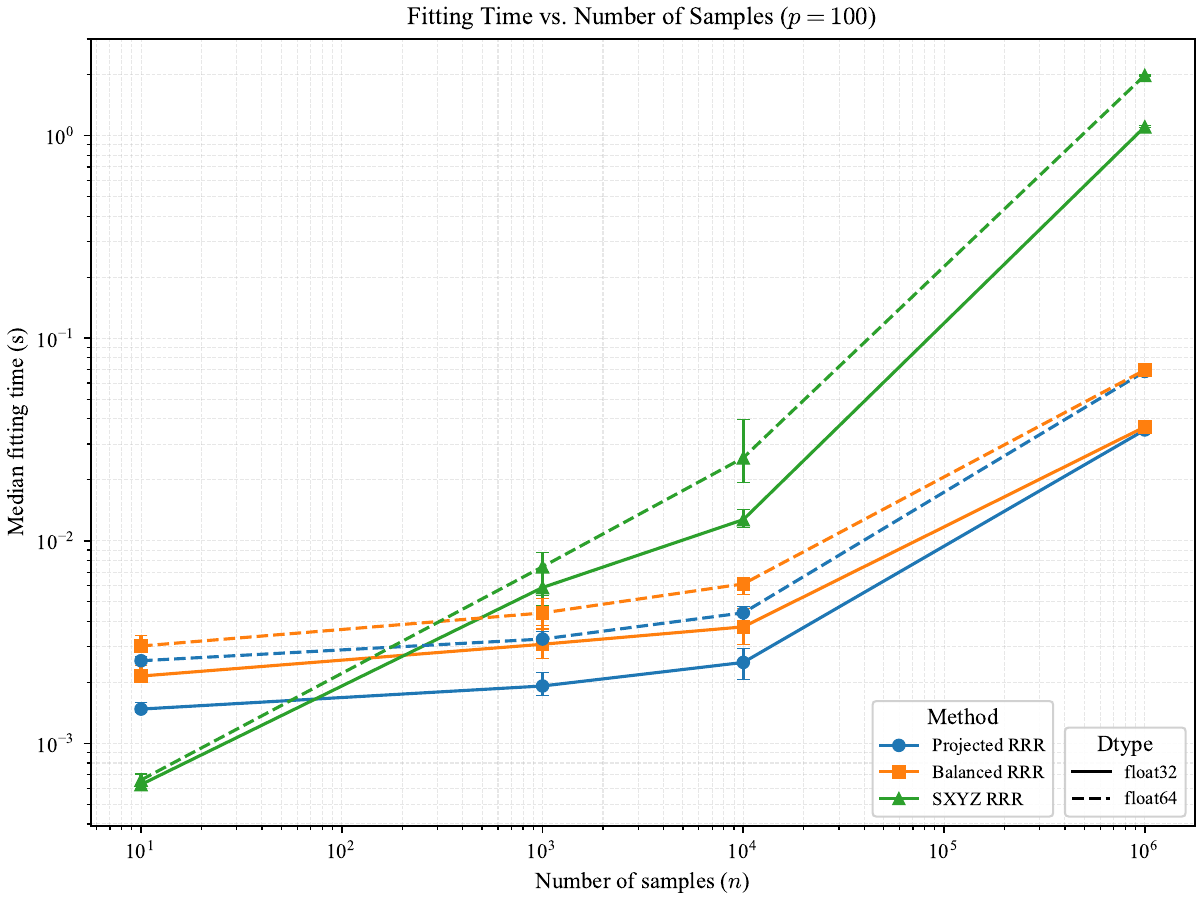}
    \end{center}
    \caption{Fitting time of the 3 algorithms on CPU for different $n$ and $p=q=100$, aggregated as in \Cref{fig:fitting_time_cpu}.}
    \label{fig:fitting_time_cpu_100}
\end{figure}

\begin{figure}
    \begin{center}
        \includegraphics[width=0.9\textwidth]{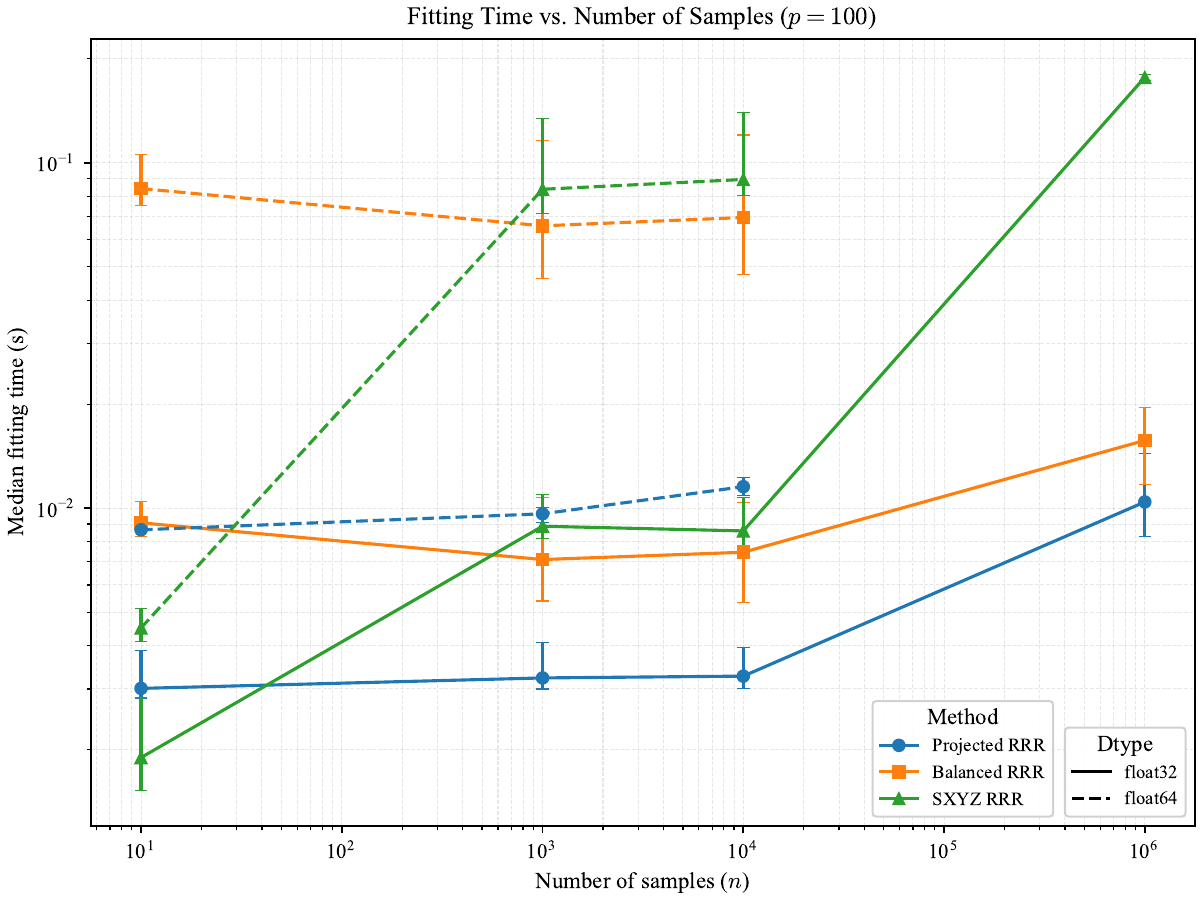}
    \end{center}
    \caption{Fitting time of the 3 algorithms on GPU for different $n$ and $p=q=100$, aggregated as in \Cref{fig:fitting_time_cpu}. The float64 runs with $n=10^6$ were not performed on GPU.}
    \label{fig:fitting_time_gpu_100}
\end{figure}

\begin{figure}
    \begin{center}
        \includegraphics[width=0.9\textwidth]{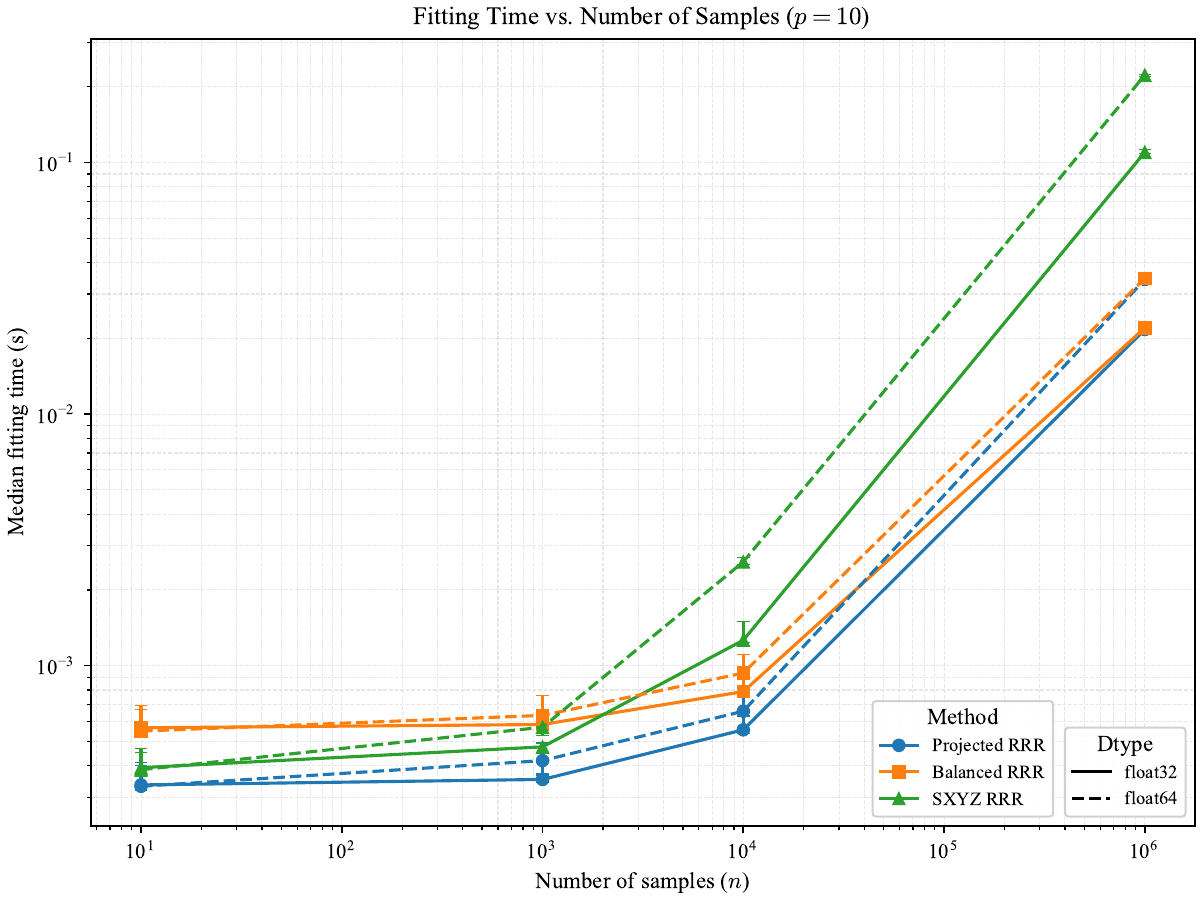}
    \end{center}
    \caption{Fitting time of the 3 algorithms on CPU for different $n$ and $p=q=10$, aggregated as in \Cref{fig:fitting_time_cpu}.}
    \label{fig:fitting_time_cpu_10}
\end{figure}

\begin{figure}
    \begin{center}
        \includegraphics[width=0.85\textwidth]{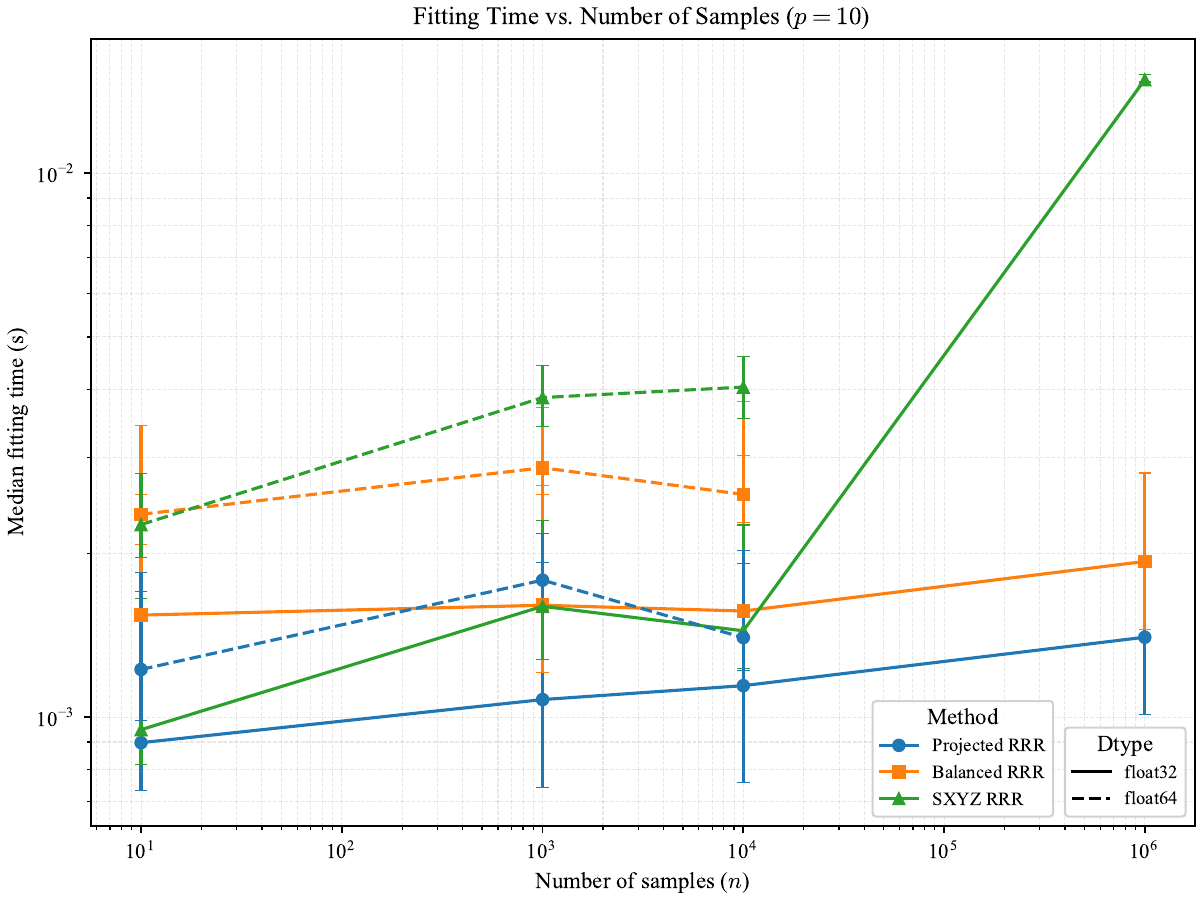}
    \end{center}
    \caption{Fitting time of the 3 algorithms on GPU for different $n$ and $p=q=10$, aggregated as in \Cref{fig:fitting_time_cpu}. The float64 runs with $n=10^6$ were not performed on GPU.}
    \label{fig:fitting_time_gpu_10}
\end{figure}

\end{document}